\documentclass[a4paper,notitlepage, 8pt,reqno]{article}

\usepackage[cp1251]{inputenc}
 \usepackage[T2A]{fontenc}
 \usepackage[russian,english]{babel}
\usepackage[all,arc,poly,2cell,curve,arrow,tips]{xy}
\usepackage{enumerate, eucal, amsthm,amsmath, amssymb}
\usepackage{graphicx}
\usepackage{mathtext}
\usepackage{amsmath,tikz}

\usetikzlibrary{tikzmark}

\usepackage{tikz}

\allowdisplaybreaks[4]

\theoremstyle{definition}

\theoremstyle{plain}
\newtheorem{thm}{Theorem}
\newtheorem*{thm*}{Theorem}

\newtheorem*{prop*}{Proposition}

\newtheorem{lem}{Lemmа}
\newtheorem{remark}{Remark}

\renewcommand{\abstractname}{}

\title{ A functional approach to a Gelfand-Tsetlin type base for   $\mathfrak{o}_5$ }

\author{D.V. Artamonov}
\date{}

\begin{document}
    \maketitle

 \maketitle

\renewcommand{\abstractname}{}

\begin{abstract}
A realization of representations of the Lie algebra $\mathfrak{o}_5$ in the space of functions on a group $Spin_5\simeq Sp_4$ is considered. In a representation we take a Gelfand-Tsetlin type base associated with a restriction $\mathfrak{o}_5\downarrow\mathfrak{o}_3$. Such a base is useful is problems appearing in quantum mechanics. We construct explicitely functions on the group that correspond to base vectors. As in the cases of Lie algebras $\mathfrak{gl}_3$, $\mathfrak{sp}_4$ these functions can be expressed through $A$-hypergeometric functions (this does not hold for algebras of these series in higher dimentions). Using this realization formulas for the action of generators are obtained.
\end{abstract}

\section{Introduction}

In the paper  \cite{g1} I.M. Gelfand and M.L.  Tsetlin  constructed a base in an irreducible finite dimensional representation of  $\mathfrak{o}_N$.  A construction of this base  is based on investigation of branching of an irreducible representation of   $\mathfrak{o}_N$ under a restriction of algebras $\mathfrak{o}_N\downarrow\mathfrak{o}_{N-1}$.  Later  it was understood that it is  reasonable t also  consider  a Gelfand-Tsetlin type base in representations  of  $\mathfrak{o}_{2n+1}$ using restrictions  $\mathfrak{o}_{2n+1}\downarrow\mathfrak{o}_{2n-1}$  i.e. restrictions  inside the series  $B$.  Such a base was constructed by Molev in  \cite{m}. He also obtained formulas for the action of generators of the algebra  $\mathfrak{o}_{2n+1}$  in this base.  Mention that construction of such a Gelfand-Tsetlin type base is not unique.

 There a lot of papers where authors construct  a Gelfand-Tsetlin type base for
$\mathfrak{o}_5$ by hand  and use it to classify states of a five-dimensional quasi-spin in the nuclear shell model  \cite{f1},\cite{f2},\cite{f3}.  Such construction give a base other than constructed in    \cite{m}.  Let us mention  also the paper   \cite{asf} that is also motivated by problems coming from quantum mechanics. In these papers  the  Gelfand-Tsetlin type base for
$\mathfrak{o}_5$ was constructed using a technique of the projection operator  (see definition in  \cite{t1}). This operator  plays an important  role in the theory of the Gelfand-Tsetlin bases (see a review  \cite{t2}),  and also in general constructions in  \cite{m}.


In the present paper we construct the Gelfand-Tsetlin type base for  $\mathfrak{o}_5$, but we use new ideas  (other than in  \cite{m}).
A starting point is the paper  \cite{bb1963}, where in the case  $\mathfrak{gl}_3$ the following observation is done. If one uses a functional realization  then a function of the group  $GL_3$, corresponding to Gelfand-Tsetlin base vectors can be expressed through a Gauss hypergeometric function  (see a modern approach in  \cite{a2}).  This fact was used in   \cite{a1}  to obtain an explicit formula for an arbitrary Clebsh-Gordan coefficient for  $\mathfrak{gl}_3$.

Explicit formulas for the Gelfand-Tsetlin base vectors for   $\mathfrak{sp}_4$ (such base vectors were constructed in  \cite{zh})  were obtained in  \cite{a4}. In this paper the formulas for the action of generators were derived.


Mention that one has an isomorphism of Lie algebras $\mathfrak{sp}_4\simeq\mathfrak{o}_5$,  thus a representation $V$ of the algebra $\mathfrak{ o}_5$ is automatically a representation of  $\mathfrak{sp}_4$ But under this isomorphism the subalgebras $\mathfrak{sp}_2\subset \mathfrak{sp}_4$ and  $\mathfrak{o}_3\subset \mathfrak{o}_5$  are not identified. Hence the Gelfand-Tsetlin type bases in  $V$ corresponding to these chains of subalgebras are different. In applications the base corresponding to the chain   $\mathfrak{o}_3\subset \mathfrak{o}_5$ is needed.

 Thus in the present paper the following problem is considered. Take a Lie group  with the Lie algebra   $\mathfrak{o}_5$.   For example one can take   $SO_5$ or $Spin_5=Sp_4$. Consider functions on this group. This space of functions  is a representation of   $\mathfrak{o}_5$, one has a canonical embedding of a arbitrary irreducible finite dimensional representation into this functional representation. It is natural to ask the following question: which  functions correspond to the Gelfand-Tsetlin base vectors? What is the action of generators onto them? We obtain answers to both question in the present paper.

Surprisingly, the case of the algebra  $\mathfrak{o}_5$ is much more difficult than the case of the algebra   $\mathfrak{sp}_4$.   The  first steps are similar to those in    \cite{a4}.   At the beginning it is more convenient to work with the space of functions on   $SO_5$. We establish a relation between restriction problems    $\mathfrak{o}_5\downarrow \mathfrak{o}_3$  and  $\mathfrak{gl}_3\downarrow \mathfrak{gl}_1$\footnote{A restriction problem $g\downarrow h$ (where $h\subset g$ are Lie algebras) is a problem of description of   $h$-highest vectors in an irreducible representation of  $g$}.  Using this relation we construct a base in the space  of $\mathfrak{o}_3$-highest vectors with a fixed   $\mathfrak{o}_3$-weight in a given  $\mathfrak{o}_5$-representation.
The function corresponding to  basic $\mathfrak{o}_3$-highest vectors are written as   $\Gamma$-series.

To obtain formulas for  functions corresponding to arbitrary vectors (which are not $\mathfrak{o}_3$-highest),  we need to apply a power of a lowering operator from  $\mathfrak{o}_3$. It is easy to check that one gets a difficult function which is not a   $\Gamma$-series. Thus the situation differs from the situation in the paper \cite{a4}.

To write the function in a simpler way we use an isomorphism  $\mathfrak{sp}_4\simeq\mathfrak{o}_5$   and realize our representation  in the space of functions on   $Spin_5=Sp_4$.  Then the functions corresponding to the  Gelfand-Tsetlin base vectors are written as   $\Gamma$-series.

But also there are similarities with the paper   \cite{a4}.   To derive formulas for the action of generators we use the Principle Lemma, formulated in Section   \ref{oslmma},  which gives a formula  for a result of  a multiplication of a $\Gamma$-series depending on minors of the matrix composed of matrix elements (considered as functions on the group) onto a minor.           In the Lemma we claim that a  lattice underlying  the   $\Gamma$-series is related to relations between minors.

The Gelfand-Tsetlin base in a representation of  $\mathfrak{o}_5$ based on restrictions  $\mathfrak{o}_5\downarrow\mathfrak{o}_3$ is not unique. The base constructed in the present paper differs from  bases constructed in papers cited above. A construction of the presnet paper establishes  a remarkable relation with the theory of   $A$-hypergeometric functions   (see  \cite{GG}).  From one hand this is a beautiful demonstration of the unity of mathematics. From the other hand that makes possible to find  explicit solution of problems that seem to have no easy solution in the traditional approach to the representation theory.  Thus using a relation with the theory of  $A$-hypergeometric functions in the case $\mathfrak{gl}_3$  one manages to obtain  a really explicit formula for an arbitrary Clebsh-Gordan coefficient for   $\mathfrak{gl}_3$  (see  \cite{a1},  and especially  \cite{a6}).

\section{ Preliminary facts}

\subsection{Algebras   $\mathfrak{sp}_4$ and $\mathfrak{o}_{5}$. A functional realization of  a representations}

Take a  complex linear space with coordinated indexed by  $\pm 2$, $\pm 1$.
Consider a Lie group  $Sp_4$,  which consists of linear transformations preserving the skew-symmetric form
$\omega=e_{-2}\wedge e_2+e_{-1}\wedge e_{1}$.

It's Lie algebra  $\mathfrak{ sp}_4$  is spanned by  $$f_{i,j}=E_{i,j}-sign(i)sign(j)E_{-j,-i}, \,\,\, i,j=\pm 2,\pm 1.$$

Take a  complex linear space with coordinated indexed by  $\pm 2$, $\pm 1$, $0$.
Consider a Lie group     $SO_5$,which consists of linear transformations preserving the  scalar product whose matrix has ones  on the  second diagonal and zeros on other places.

It's Lie algebra   $\mathfrak{o}_5$   is spanned by $$F_{i,j}=E_{i,j}-E_{-j,-i}, \,\,\, i,j=\pm 2,\pm 1,0.$$

On the space of complex-valued functions on the group  $G=Sp_4,SO_5$ the group itself acts by  right shifts.  Thus the space of functions is a representation of the group   $G$ and hence of it's Lie algebra.  Let  $a_i^j$  be a function of  a matrix element where   $i $ is a column index and   $j$  is a row index.  These indices are taken from   $\{\pm 2, \pm 1, 0\}$.  Compose determinants of size $1$ and  $2$:

\begin{equation}
\label{aij}
a_i=det(a_i^{-2}),\,\,\, a_{i_1,i_2}=det(a^{j}_i)_{i=i_1,i_2}^{j=-2,-1}
\end{equation}

They are also function on the group.

The Lie algebra acts onto  \eqref{aij}  by the ruler

\begin{align*}
&f_{i,j}a_{i_1}=a_{i_1\mid_{j\mapsto i }}- sign(i)sign(j)a_{i_1\mid_{-i\mapsto -j }},\\
&F_{i,j}a_{i_1}=a_{i_1\mid_{j\mapsto i }}- a_{i_1\mid_{-i\mapsto -j }},
\end{align*}

where  $i_1\mid_{j\mapsto i }$ is a substitution of   $j$ instead of   $i$,  and in the case  $j\neq i_1$ one gets  $0$. Analogously  one defines an action onto  $a_{i_1,i_2}$.  And onto a product of determinants the action is defined by the Leibnitz ruler.

One can easily prove that

\begin{equation}
  \label{stv}
a_{-2}^{m_{-2}-m_{-1}}a_{-2,-1}^{  m_{-1} }
\end{equation}

is a highest vector of weight $[m_{-2},m_{-1}]$.

 The functional realization has the following advantage: one can  write conditions that define polynomials in  \eqref{aij} that form a representation with the highest vector   \eqref{stv} (in terms of the indicator system, see \cite{zh}).



\subsection{  A mapping  $Sp_4\rightarrow SO(5)$}

The exists a covering  $Sp_4\rightarrow SO_5$ than induces an isomorphism  $\mathfrak{sp}_4 =\mathfrak{o}_5$.  Note that this isomorphism can be derived from the coincidence of root systems  $C_2$ and  $B_2$.

The mapping is constructed as follows.  Consider the space  $V=\mathbb{C}<e_{-2},e_{-1},e_{1},e_{2}>$  of the standard representation of   $Sp_4$.  Then the action of   $Sp_4$ onto  $V\wedge V$ has a five-dimensional invariant subspace  $$\mathbb{C}<v_{-2}= e_{-2}\wedge e_{-1},
 v_{1}=  e_{-2}\wedge e_{1},   v_{-1}= e_{2}\wedge e_{-1}, v_{2}=e_{1}\wedge e_{2}, v_{0}=\frac{1}{\sqrt{2}}(e_{-2}\wedge e_{2}+e_{-1}\wedge e_{1})>.$$  The action of  $Sp_4$ on this subspace has an invariant scalar product such that  $<v_{-i},v_{i}>=1$,   and other scalar products vanish. Thus  the elements of  $Sp_4$ act on this space by orthogonal matrices   $5\times 5$, hence one has  $Sp_4\rightarrow SO_5$.   This mapping is an epimorphism, a preimage of each elements of   $SO_5$ consists of two elements  (in particular a preimage of   $E\in SO_5$  consists of $E$   and   $-E$ from  $Sp_4$).
Thus one has a mapping  $Fun(SO_5)\rightarrow Fun(Sp_4)$.

The obtained mapping between Lie algebras look as follows.
 The positive root elements of  $\mathfrak{sp}_4$ are mapped by the following ruler

   \begin{align*}
 \begin{array}{|c|c|c|c|}
 \hline
 \frac{1}{\sqrt{2}}f_{-2,-1} &\frac{1}{2} f_{-2,2} & - \frac{1}{\sqrt{2}} f_{-2,1} & \frac{1}{2} f_{-1,1} \\
 \hline
F_{-1,0} &  F_{-2,1} & F_{-2,0} & F_{-2,-1}\\
 \hline
 \end{array}
 \end{align*}

The negative  root elements of  $\mathfrak{sp}_4$  are mapped to corresponding negative root elements of  $\mathfrak{o}_5$.

 For the  Cartan elements one has

    \begin{align*}
 \begin{array}{|c|c|}
 \hline
 \frac{1}{2}(f_{-2,-2}+f_{-1,-1} )& \frac{1}{2}(f_{-2,-2}-f_{-1,-1}) \\
 \hline
 F_{-2,-2} &  F_{-1,-1} \\
 \hline
 \end{array}
 \end{align*}

In particular a vector  with a  $\mathfrak{sp}_4$-weight $[x,y]$ has a $\mathfrak{o}_5$-weight  $\frac{1}{2}[x+y,x-y]$.

 From an explicit form of the mapping we conclude that  $a_i$ are transformed to    $b_{j_1,j_2}$,  more precise one has:

  \begin{align}
  \label{t1}
 \begin{array}{|c|c|c|c|c|}
 \hline
 a_{-2} & a_{-1} & a_0 & a_{1}& a_{2} \\
 \hline
 b_{-2,-1} & b_{-2,1} &\frac{1}{\sqrt{2}}( b_{1,-1}+b_{-2,2}) & b_{2,-1} & b_{1,2}\\
 \hline
 \end{array}
 \end{align}

 Now let us find images of  determinants    $a_{i_1,i_2}\in Fun(SO_5)$.   These determinants form a representation with the highest weight  $[1,1]$.  Hence images of these functions in   $Fun(Sp_4)$ form a representation with the highest weight   $[2,0]$,   which is spanned by  $b_{j_1}b_{j_2}$.

 The highest vectors are  $a_{-2,-1}$  and $(b_{-2})^2$. Let us identify them then using the identification of   $\mathfrak{o}_5$ and   $\mathfrak{sp}_4$ one obtains a correspondence

 \begin{tiny}
   \begin{align}
   \label{t2}
 \begin{array}{|c|c|c|c|c|c|c|c|c|c|}
 \hline
  (b_{-2})^2 & \sqrt{2} b_{-2}b_{-1} &  \sqrt{2}b_{-2}b_1 & b_{-2}b_{2}+b_{-1}b_{1}&
b_{-2}b_{2}-b_{-1}b_{1} &\sqrt{2} b_{-1}b_{2} & -\sqrt{2} b_{1}b_{2} & (b_{-1})^2  &- (b_1)^2 & (b_2)^2\\
 \hline
 a_{-2,-1} & a_{-2,0} & a_{-1,0} & a_{-2,2} & a_{-1,1}
 &a_{0,1}&a_{0,2}&a_{-2,1}&a_{2,-1}&a_{1,2}\\
 \hline
 \end{array}
 \end{align}
\end{tiny}

Obviously the correspondence respects products.

\subsection{  $A$-hypergeometric functions}

Information about a $\Gamma$-series can be found in \cite{GG}.

Let $\mathcal{B}\subset \mathbb{Z}^N$ be a lattice  and let  $\mu\in \mathbb{Z}^N$ be a fixed vector. We call it {\it the shift vector}. Define a  {\it  hypergeometric
	$\Gamma$-series }  in variables $z_1,...,z_N$ by the formula

\begin{equation}
\label{gmr}
\mathcal{F}_{\gamma}(z,\mathcal{B})=\sum_{v\in
	\mathcal{B}}\frac{z^{v+\gamma}}{\Gamma(v+\gamma+1)},
\end{equation}
where  $z=(z_1,...,z_N)$. We use a multi-index notation:

$$
z^{v+\gamma}:=\prod_{i=1}^N
z_i^{v_i+\gamma_i},\,\,\,\Gamma(v+\gamma+1):=\prod_{i=1}^N\Gamma(v_i+\gamma_i+1).
$$

Note that  when at least one components of $v+\gamma$ is integer and negative that the coresponding summand in  \eqref{gmr} vanishes.    This is the reason, why the  considered in the present paper    $\Gamma$-series are finite sums. For simplicity we write factorials instead of  $\Gamma$-functions.

We need the following ruler for differentiation of a   $\Gamma$-series:

\begin{equation}
\frac{\partial}{\partial z_i}\mathcal{F}_{\gamma}=\mathcal{F}_{\gamma-e_i},
\end{equation}

where  $e_i$  is a standard base vector for the   $i$-th coordinate.

An $A$-hypergeometric function satisfies the Gelfand-Kapranov-Zelevinsky system of PDE, which consists of equations of two types .

{\bf 1.} Let   $\alpha=(\alpha_1,...,\alpha_N)$ be a vector orthogonal to  $\mathcal{B}$, then

\begin{equation}
\label{e1}
\alpha_1z_1\frac{\partial}{\partial z_1}\mathcal{F}_{\gamma}+...+\alpha_Nz_N\frac{\partial}{\partial z_N}\mathcal{F}_{\gamma}=(\alpha_1\gamma_1+...+\alpha_N\gamma_N)\mathcal{F}_{\gamma},
\end{equation}
it is sufficient to consider base vectors of the orthogonal compliment to   $\mathcal{B}$.

{\bf 2.}  Ler  $v\in \mathcal{B}$ and  $v=v_+-v_-$,  where vectors  $v_+$, $v_-$ have non-negative coordinates.  Take non-zero elements
$v_+=(...v_{i_1},....,v_{i_k}...)$,  $v_-=(...v_{j_1},....,v_{j_l}...)$. Then

\begin{equation}
\label{e2} (\frac{\partial }{\partial
	z_{i_1}})^{v_{i_1}}...(\frac{\partial}{\partial z_{i_k}})^{v_{i_k}}
\mathcal{F}_{\gamma}=(\frac{\partial }{\partial
	z_{j_1}})^{v_{j_1}}...(\frac{\partial }{\partial z_{j_l}})^{v_{j_l}} \mathcal{F}_{\gamma}.
\end{equation}

It is sufficient to consider only base vectors  $v\in \mathcal{B}$.

\subsection{A Gelfand-Tsetlin type base for a chain  $g\supset  h$}
\label{bgc}

Let us introduce a unified notation for Lie algebras: denote as   $g$  the algebra  $\mathfrak{o}_5$ of   $\mathfrak{sp}_4$,  and denote as  $h$ the subalgebra  $<F_{0,-2},F_{-2,-2},F_{-2,0}>$ or $<f_{-2,1},f_{-2,-2}-f_{1,1},f_{1,-2}>$.

Let   $V$  be a standard representation of   $g=\mathfrak{o}_5$ with the highest weight  $\mu_{2}=[m_{-2},m_{-1}]$.  Take a subalgebra $h=\mathfrak{o}_3=<F_{0,-2},F_{-2,-2},F_{-2,0}>$ and consider   $V$ as a representation of this algebra.     Now it is not irreducible but it splits into a sum of irreducibles. Consider the space  of   $\mathfrak{o}_3$-highest vectors  with the highest weight $\mu_{1}=[s_{-2}]$.  Chose a base in this space and let  $\mu'_{2}$ be an index numerating the base vectors. In the space of a  $\mathfrak{o}_3$-representation generated by  a chosen $\mathfrak{o}_3$-highest vector  chose a weight  base  and let  $\mu'_1$ be an index numerating the base vectors. Then one gets that in the space     $V$ there exists a base indexed by

\begin{equation}
\label{dc}
\begin{pmatrix}
\mu_2\\
\mu'_2\\
\mu_1\\
\mu'_1
\end{pmatrix}
\end{equation}

Such a base is called a Gelfand-Tsetlin type  base. Since the choice of a base in the space of  $\mathfrak{o}_3$-highest vectors is not unique, the construction of such a base in not unique.

One can describe  this construction on the language of   $ g=\mathfrak{sp}_4$. In this case one must consider a subalgebra  $h=<f_{-2,1},f_{-2,-2}-f_{1,1},f_{1,-2}>$,  isomorphic to  $\mathfrak{sl}_2$.


\section{ A base in the space of   $h$-highest vectors of a representation of  $g$}

Let us be given a representation   $V$ of the Lie algebra   $g$. The first step in the construction of the Gelfand-Tsetlin type  base in $V$  for the chain $g\supset h$ is a construction of base in the space of   $h$-highest vectors in $g$.

It is more natural to do it on the language of  $\mathfrak{o}_5$, however to finish the construction of the Gelfand-Tsetlin type  base   we need to pass to the language of   $\mathfrak{sp}_4$.

\subsection{A base in the space of    $\mathfrak{o}_3$-highest vectors }
\label{o335}
In   \cite{zh}  an explicit description of conditions is given which define in the space of functions an irreducible representation with the highest vector    \eqref{stv}.    They look as follows

\begin{enumerate}
\item $L_{-}f=0$, where  $L_{-}$ is a left infinitesimal shift  by an element of    $SO_{5}$, corresponding to a negative root .
\item $L_{-i,-i}f=m_{-i}f$, where   $L_{-i,-i}$, $i=1,2$  is a left infinitesimal shift  by an element of
$SO_{5}$, corresponding to a Cartan element  $F_{-i,-i}$.
\item $f$ satisfies the indicator system 

$$\begin{cases}L_{-2,-1}^{l_{-2}+1}f=0,\,\,\, l_{-2}=m_{-2}-m_{-1}    , \\L_{-1,0}^{l_{-1}+1}f=0\,\,\, l_{-1}=2m_{-1}\end{cases},$$
 where   $L_{i,j}$  is a left infinitesimal shift  by an element  $F_{i,j}$.
\end{enumerate}

Consider the space of   $\mathfrak{o}_3$-highest vectors with a given   $\mathfrak{o}_3$-weight.  A function is a   $\mathfrak{o}_3$-highest vector if it depends on  \begin{equation}\label{a3}a_{-2},a_{-1},a_{1},a_{-2,-1},a_{-2,1},a_{-1,1}.\end{equation}

Let us find among function depending on these variables the functions satisfying  conditions  1-3.   Let us introduce for shots time a new notation for determinants which includes upper indices of columns participating in the determinant (see  \eqref{aij})

$$
a^{-2}_{i},\,\,\,\,a^{-2,-1}_{i_1,i_2}.
$$
One can easily prove that   $L_{-i,-j}$ acts as follows:

\begin{align}
\begin{split}
\label{lij}
&L_{-2,-2}a^{-2}_{i_1}=a^{-2}_{i_1},\,\,\,L_{-1,-1}a^{-2}_{i_1}=0\,\,\,\,L_{-i,-i}a^{-2,-1}_{i_1,i_2}=a^{-2,-1}_{i_1,i_2},\\
&L_{-2,-1}a_{i_1}^{-2}=a_{i_1}^{-1},\,\,\,L_{-2,-1}a_{i_1}^{-1}=0,\,\,\,L_{-1,0}a_{i_1}^{-2}=0,\\
&L_{-2,-1}a_{i_1,i_2}^{-2,-1}=0,\,\,\,L_{-1,0}a_{i_1,i_2}^{-2,-1}=a_{i_1,i_2}^{-2,0},\,\,\,L_{-1,0}a_{i_1,i_2}^{-2,0}=a_{i_1,i_2}^{-2,1},\,\,\,
L_{-1,0}a_{i_1,i_2}^{-2,1}=0
\end{split}
\end{align}

Among functions depending on   \eqref{a3},  let us find functions satisfying the conditions 1-3.

\begin{thm}
\label{lempoc}   A function depending on  \eqref{a3} is a  $\mathfrak{o}_3$-highest vector in an irreducible representation with the highest vector \eqref{stv}   if and only if the following condition holds.

If the highest weight is integer then $f$ is a polynomial in determinants and the sum of exponents of determinants   of order  $i$ equals $l_{-i}$.

If the highest weight is half-integer then \begin{equation}\label{f12}f=(a_{-2,-1})^{\frac{1}{2}}f_1+a_{-2,0}(a_{-2,-1})^{-\frac{1}{2}}f_2,\end{equation} where  $f_1$ and  $f_2$ are polynomials in determinants.  The sum of exponents of determinants   of order $1$ in each summand of   $f$  equals  $l_{-2}$,   and the sum of exponents of determinants   of order $2$ equals $l_{-1}$.

\end{thm}

\proof

Consider first the case of an integer highest weight, the vector  \eqref{stv}  is a polynomial in determinants and under the action of elements of the  algebra  the space of polynomials is invariant. Thus an arbitrary vector of the representation is written as a polynomial in determinants. Then the statement of the Theorem follows from  \eqref{lij}.

Now consider the case of the half-integer highest weight.  First of all let us show that a functions satisfying the conditions 1-3, 
looks as follows
\eqref{f12}.

Consider the case of the highest vector  $(a_{-2,-1})^{\frac{1}{2}}$.  An arbitrary vector of the representation is a linear combination of vectors

$$
F_{0,-1}^{p_{-1}}F_{-1,-2}^{p_{-2}}(a_{-2,-1})^{\frac{1}{2}}.
$$

Such a vector is nonzero only if  $p_{-2}=0$,  and  $p_{-1}=0$
or  $1$. Indeed if   $p_{-1}=1$ one gets the vector
$\frac{1}{2}a_{-2,0}(a_{-2,-1})^{-\frac{1}{2}}$, and for
$p_{-1}=2$ one gets

\begin{align*}
&-\frac{1}{2}a_{-2,1}(a_{-2,-1})^{-\frac{1}{2}}-\frac{1}{4}(a_{-2,0})^2(a_{-2,-1})^{-\frac{3}{2}}=\\
&=-\frac{1}{2}(a_{-2,-1})^{-\frac{3}{2}}(a_{-2,1}a_{-2,-1}+\frac{1}{2}(a_{-2,0})^2)=0,
\end{align*}

we have used a relations  $a_{-2,0}^2=-\frac{1}{2}a_{-2,1}a_{-2,-1}$.

Now consider the case of a half-integer highest weight. The highest vector  \eqref{stv}
can be written as follows:

$$
v_0=v'_0(a_{-2,-1})^{\frac{1}{2}},
$$

where  $v'_0$  is a polynomial in determinants. An arbitrary vector $f$
of the representation is obtained as a linear combination of vectors  that are results of application  to the highest vectors of the operators   $
F_{0,-1}^{p_{-1}}F_{-1,-2}^{p_{-2}}$.  As a result one gets a vector of type
\eqref{f12}.

Now we need to show that each vector of type \eqref{f12},
satisfying the conditions from the formulation of the Theorem satisfies conditions  1-3.
Note that conditions 1,2 are satisfied automatically one needs to check only the condition  3.

The operator $L_{-2,-1}$  acts onto determinants of order $1$ only. Such determinants are included only in   $f_1$, $f_2$,  thus they occur in non-negative integer powers, the sum of powers of determinants of  order  $1$ equals $l_{-2}$. Hence the condition
$L_{-2,-1}^{l_{-2}+1}f=0 $  holds.

Now consider the equation
$L_{-1,0}^{2m_{-1}+1}f=L_{-1,0}^{2[m_{-1}]+1+1}f=0$, where
$[m_{-1}]$ is an integer part.  The operator $L_{-1,0}^{2[m_{-1}]+1+1}$ acts onto each summand  in \eqref{f12} by the Leibnitz ruler as follows.

 Either $L_{-1,0}^{2[m_{-1}]+1+1}$ acts onto the second factor   $f_1$  or  $f_2$ only. Then one obtains   $0$,  since in  $f_1$ and  $f_2$  the sum of powers of determinants of order   $2$ equals  $[m_{-1}]$,  but such polynomials are annihilated by  $L_{-1,0}^{2[m_{-1}]+1}$.

 Either $L_{-1,0}^{2[m_{-1}]+1}$ acts onto the second factor   $f_1$ or  $f_2$,  and
   $L_{-1,0}$ acts onto the first factor. Then one gets   $0$  by the same reason.

 Either $L_{-1,0}^{2[m_{-1}]+2-k}$ acts onto the second factor,  and  the operator  $L_{-1,0}^{k}$, where $k\geq 2$ acts onto the first factor. Since the first  factor is a vector of representation with the highest weight $[\frac{1}{2},\frac{1}{2}]$, then under the action  $L_{-1,0}^{2}$
it vanishes.

 Hence the vector of type \eqref{f12} vanishes under the action of
$L_{-1,0}^{2m_{-1}+1}$.

\endproof

From the results of  \cite{m} it follows that in the space  of  $\mathfrak{o}_3$-highest vectors of a
$\mathfrak{o}_5$-representation with a fixed  $\mathfrak{o}_3$-weight there is a base indexed by the number $\sigma$
  and an integer or half-integer diagram satisfying the betweeness conditions

\begin{align}
\begin{split}
\label{dio5}
&\sigma,\begin{pmatrix}m_{-2}&& m_{-1} && 0\\
 &k_{-2}&&   k_{-1}  \\
&&s_{-2}
\end{pmatrix},
\end{split}
\end{align}
where  $\sigma=0,1$. In the case $k_{-1}=0$, one has   $\sigma=0$.

Below in the same space  we construct  another base whose elements are indexed by the same diagrams.  For this new base we manage to write explicitly the functions corresponding to diagrams.

Consider the space  $\mathbb{C}^4$.  Let us write the determinants in the following order:  $$a=(a_{-2},a_{-1},a_{-2,1},a_{-1,1}).$$ To the diagram \eqref{dio5} the corresponds the vector  $$\gamma=(s_{-2}-m_{-1},k_{-2}-s_{-2},m_{-1}-k_{-1},0).$$ Introduce a lattice
$$\mathcal{B}^1=\mathbb{Z}<(1,-1,-1,1)>.$$

\begin{lem}
	\label{lm1}
	Let us be given a representation of   $\mathfrak{o}_5$ with the highest weight  $[m_{-2},m_{-1}]$. In the space of    $\mathfrak{o}_3$-highest vectors with weight  $s_{-2}$ there exists a base indexed by diagrams of type  \eqref{dio5}. In this diagram in the case of the integer highest weight if   $\sigma=1$ then  $k_{-1}\geq 1$.
To the diagram \eqref{dio5} there corresponds the function \begin{equation}\label{bbo}a_{-2,0}^{\sigma}
a_{1}^{m_{-2}-k_{-2}}a_{-2,-1}^{k_{-1}-\sigma}\mathcal{F}_{\gamma}(a,\mathcal{B}^1).\end{equation}

The weights are calculated as follows: $s_{-2}$ is a 
$(-2)$-component of the weight, a
   $(-1)$-component of the weight equals  $-2(k_{-2}+k_{-1})+(m_{-2}+m_{-1})+s_{-2}+\sigma$.

\end{lem}

\begin{remark}

The function is written explicitly as follows
\begin{equation}
a_{-2,0}^{\sigma}a_{1}^{m_{-2}-k_{-2}}a_{-2,-1}^{k_{-1}-\sigma}\sum
\frac{1}{p_{-1}!p_{-1,1}!p_{1}!p_{-2,1}!} a_{-1}^{p_{-1
}}a_{-1,1}^{p_{-1,1}}a_{-2}^{p_{-2}}a_{-2,1}^{p_{-2,1}},
\end{equation}
the summation is taken over all integer non-negative collections  $p_{-1},p_{-1,1},p_{1},p_{-2,0}$, such that

 \begin{align*}
&p_{-2}+p_{-1}=k_{-2}-m_{-2},\,\,\,
p_{-1}+p_{-2}=k_{-2}-m_{-1},\,\,\,\\&p_{-1,1}+p_{-2,1}=m_{-1}-k_{-1},\,\,\, p_{-1}+p_{-1,1}=k_{-2}-s_{-2}.
\end{align*}


 \end{remark}
\proof
Let us prove that such functions form a base in the considered linear space.

Consider first the case of the integer highest weight. The functions described in Theorem \ref{lempoc} are polynomials in determinants.  Consider the admissible polynomial  $f$ in determinants satisfying the conditions of Theorem  \ref{lempoc}.
  Then  $f$ as an element of the representation with the highest vector  $v_0$ can be written as a linear combination of polynomials of type   $F_{-1,-2}^{p}F_{0,-1}^{q}v_0$, where $v_0=a_{-2}^{m_{-2}-m_{-1}}a_{-2,-1}^{m_{-1}}$.
One has \begin{equation}\label{f01}F_{0,-1}^2a_{-2,-1}=2F_{0,-1}a_{-2,-1}a_{-2,0}=2a_{-2,0}^2-2a_{-2,-1}a_{-2,1}=-4a_{-2,-1}a_{-2,1},\end{equation} using it in the case when   $q=2q'$ is even one gets 

$$F_{-1,-2}^{p}F_{0,-1}^{2q'}v_0=const E_{-1,-2}^pE_{1,-1}^{q'}v_0.$$

We formally define an action of  $E_{i,j}$ onto symbols  $a_{i_1}$, $a_{i_1,i_2}$ by the ruler  $$E_{i,j}a_X=a_{X\mid_{j\mapsto i}}.$$

The polynomial on the right satisfies the conditions of Theorem   \ref{lempoc} for the algebra  $\mathfrak{gl}_3$ and the highest weight
$[m_{-2},m_{-1},0]$. Thus one obtains an isomorphism between the  span of vectors $F_{-1,-2}^{p}F_{0,-1}^{2q'}v_0$ the space of representation of  $\mathfrak{gl}_3$ with the highest weight
$[m_{-2},m_{-1},0]$.

Let us use the Biedenharn-Baird Theorem  (see  \cite{bb1963}, and it's modern formulation in  \cite{a2}).  It states that  in the space of representation $\mathfrak{gl}_3$ with the highest weight 
$[m_{-2},m_{-1},0]$ in the functional realization the Gelfand-Tsetlin vector corresponding to the diagram   \eqref{dio5} is  written as  the function $\mathcal{F}_{\gamma}(a,\mathcal{B}^1)$.

Thus in the span of vectors $F_{-1,-2}^{p}F_{0,-1}^{2q'}v_0$ there exists a base \eqref{bbo},
given by diagrams \eqref{dio5}, where $\sigma=0$.

Note that in the considered case the eigenvalues of  $F_{-2,-2}$ and   $F_{-1,-1}$  correspond to eigenvalues  of  $E_{-2,-2}$ and $E_{-1,-1}-E_{1,1}$. They are equal to   $s_{-2}$ and $-2(k_{-2}+k_{-1})+(m_{-2}+m_{-1})+s_{-2}$.

In the case when  $q=2q'+1$ is odd one has

 $$F_{-1,-2}^{p}F_{0,-1}^{2q'+1}(a_{-2}^{m_{-2}-m_{-1}}a_{-2,-1}^{m_{-1}})=const (E_{-1,-2}^pE_{1,-1}^{q'}(a_{-2}^{m_{-2}-m_{-1}}a_{-2,-1}^{m_{-1}-1}))a_{-2,0}.$$

If one removes  $a_{-2,0}$, than the polynomial in determinants on the right satisfies the conditions of Theorem    \ref{lempoc} for  $\mathfrak{gl}_3$ and the highest weight $[m_{-2}-1,m_{-1}-1,0]$.
 In the space of such polynomials there exists the Gelfand-Tsetlin base indexed by integer diagrams

 \begin{align}
 \begin{split}
\label{ggoc3}
\begin{pmatrix}m_{-2}-1 &&  m_{-1}-1 && 0\\
&k_{-2}-1& & k_{-1}-1  \\
&&s_{-2}-1\end{pmatrix}.
\end{split}
\end{align}
Note that here $k_{-1}-1\geq 0$.

In the span  $F_{-1,-2}^{p}F_{0,-1}^{2q'+1}v_0$ there exists a base \eqref{bbo}, given by diagrams \eqref{dio5}, where $\sigma=1$ in the case $k_{-1}\geq 1$.

In the odd case the eigenvalues of $F_{-2,-2}$ and $F_{-1,-1}$ correspond to eigenvalues of   $E_{-2,-2}-1$ and $E_{-1,-1}-E_{1,1}$. They are equal to  $s_{-2}-1$ and $-2(k_{-2}-1+k_{-1}-1)+(m_{-2}-1+m_{-1}-1)+s_{-2}-1=-2(k_{-2}+k_{-1})+(m_{-2}+m_{-1})+s_{-2}+1$.


Thus in the case of integer highest weight there exist a base   \eqref{bbo}, indexed by \eqref{dio5}  where  $\sigma=0,1$.

Now consider the case of a half-integer highest weight.  One has
\begin{equation}\label{e38}F_{0,-1}a_{-2,-1}^{1/2}=\frac{1}{2}a_{-2,0}a_{-2,-1}^{-1/2},\,\,\,\, F_{0,-1}^{2}a_{-2,-1}^{1/2}=0.\end{equation}
The highest vector looks as follows $v_0=a_{-2}^{m_{-2}-m_{-1}}a_{-2,-1}^{[m_{-1}]+1/2}$, where  $[m_{-1}]$ is an integer part.  A vector of the representation can be written as a linear combination of vectors of type

$$F_{-1,-2}^{p}F_{0,-1}^{q}(a_{-2}^{m_{-2}-m_{-1}}a_{-2,-1}^{[m_{-1}]+1/2}).$$

If  $q=2q'$  is even then using \eqref{f01}, \eqref{e38}, one gets that  the previous expression can be written as:
 \begin{align*}
&F_{-1,-2}^{p}F_{0,-1}^{2q'}(a_{-2}^{m_{-2}-m_{-1}}a_{-2,-1}^{[m_{-1}]+1/2})=const E_{-1,-2}^{p}E_{1,-1}^{q'}(a_{-2}^{m_{-2}-m_{-1}}a_{-2,-1}^{[m_{-1}]})a_{-2,-1}^{1/2}
.\end{align*}

Thus there exists a natural isomorphism  between the span of vectors of type $F_{-1,-2}^{p}F_{0,-1}^{2q'}v_0$ and the space of vectors in the representation  $\mathfrak{gl}_3$ with the highest weight
$[m_{-2}-\frac{1}{2},m_{-1}-\frac{1}{2},0]$.  In this space there exists a base indexed by integer diagrams of type

\begin{align}
 \begin{split}
\label{gggoc3}
\begin{pmatrix} m_{-2}-\frac{1}{2} && m_{-1}-\frac{1}{2} && 0\\
&k_{-2}-\frac{1}{2}&&  k_{-1}-\frac{1}{2}  \\
&&s_{-2}-\frac{1}{2}
\end{pmatrix}.
\end{split}
\end{align}

Using the formula for the vector corresponding to the diagram in the case $\mathfrak{gl}_3$, one gets that in the span
$F_{-1,-2}^{p}F_{0,-1}^{2q'}v_0$ there exists a base  \eqref{bbo}, indexed by half-integer diagrams of type \eqref{dio5} with $\sigma=0$.

The eigenvalues of $F_{-2,-2}$ and $F_{-1,-1}$
 correspond to eigenvalues of  $E_{-2,-2}-\frac{1}{2}$ and
$E_{-1,-1}-E_{1,1}-\frac{1}{2}$. One easily proves that $s_{-2}$
is a   $(-2)$-component of the weight, and a 
$(-1)$-component of the weight equals
$-2(k_{-2}+k_{-1})+(m_{-2}+m_{-1})+s_{-2}$.

Now let $q=2q'+1$ be odd, one has

\begin{align*}
&F_{-1,-2}^{p}F_{0,-1}^{2q'+1}
(a_{-2}^{m_{-2}-m_{-1}}a_{-2,-1}^{[m_{-1}]+1/2})=const (
E_{-1,-2}^{p}E_{1,-1}^{q'}
a_{-2}^{m_{-2}-m_{-1}}a_{-2,-1}^{[m_{-1}]})a_{-2,0}a_{-2,-1}^{-1/2}
,\end{align*}

thus one gets an isomorphism between the span vectors 
$F_{-1,-2}^{p}F_{0,-1}^{2q'+1}v_0$ and the space of the representation of 
$\mathfrak{gl}_3$ with the highest weight
$[m_{-2}-\frac{1}{2},m_{-1}-\frac{1}{2},0]$.  Thus in the span there exists the base \eqref{bbo}, indexed by integer diagrams of type \eqref{dio5} with  $\sigma=1$. The eigenvalues of 
$F_{-2,-2}$ and $F_{-1,-1}$ correspond to eigenvalues of 
$E_{-2,-2}-\frac{1}{2}$ and  $E_{-1,-1}-E_{1,1}+\frac{1}{2}$. One easily proves that  $s_{-2}$ is a   $(-2)$-component of the weight, and a
$(-1)$-component of the weight equals
$-2(k_{-2}+k_{-1})+(m_{-2}+m_{-1})+s_{-2}+1$.

Thus in the case of the half-integer highest weight these exists a base   \eqref{bbo}, indexed by diagrams  \eqref{dio5}, where $\sigma=0,1$.

Thus we have proved that in the space of functions defined by Theorem  \ref{lempoc} there exists a base indexed by diagrams   \eqref{dio5}.  Such functions define  $\mathfrak{o}_3$-highest vectors in a given irreducible representation of  $\mathfrak{o}_5$. According to   \cite{m}, in this space there exist a base indexed by  \eqref{dio5}. Thus using the arguments of dimension   we obtain that we have constructed a base in the space of all  $\mathfrak{o}_3$-highest vectors.

Let us formulate the ruler of calculation of the weight: $s_{-2}$ is a 
$(-2)$-component of the weight, 
 a  $(-1)$-component of the weight  equals $-2(k_{-2}+k_{-1})+(m_{-2}+m_{-1})+s_{-2}+\sigma$.

\endproof

	








\subsubsection{Further considerations}

An arbitrary vector is obtained form a   $\mathfrak{o}_3$-highest vector corresponding to a diagram   \eqref{dio5} by application of $\frac{F_{0,-2}^p}{p!}$, $p=0,...,2s_{-2}$.  Such a vector is encoded by a diagram


\begin{equation}
\label{dc1}
\sigma,\begin{pmatrix}
m_{-2} && m_{-1} && 0\\
&k_{-2} && k_{-1}\\
&&s_{-2} && -s_{-2}\\
&&&s_{-1}
\end{pmatrix}
\end{equation}

But when one applies  $\frac{F_{0,-2}^p}{p!}$  to  $a_{-2,0}^{\sigma}a_{1}^{m_{-2}-k_{-2}}a_{-2,-1}^{k_{-2}-\sigma}\mathcal{F}_{\gamma}(a,\mathcal{B}_1)$ one gets not a   $\Gamma$-series but a much more difficult functions.

It turns out that on the language of  $\mathfrak{sp}_4$ one can construct the Gelfand-Tsetlin type base encoded by   \eqref{dc1},  such that each vector is presented by a function on  $Sp_4$, written as a  $\Gamma$-series.

\subsubsection{A base in the space of  $\mathfrak{h}$-highest vectors}
Let us write the results obtained above on the language of  $g=\mathfrak{sp}_4$ and it's subalgebra   $h$.
Let us find an image of   \eqref{bbo} in   $Fun(Sp_4)$. The mapping
 $Fun(SO_5)\rightarrow Fun(Sp_4)$ is described explicitly in   \eqref{t1}, \eqref{t2}.  One gets that the image of  the function  \eqref{bbo} up to multiplication by a constant is written as

 \begin{equation}
 \label{bfg}
 (b_{-2}b_{-1})^{\sigma}b_{-1,2}^{2(m_{-2}-k_{-2})}(b_{-2})^{k_{-1}}\cdot \mathcal{F}_{\gamma}(b_{-2,-1},b_{-2,1},(b_{-1})^2,-b_{1}b_{-1}+b_{-2}b_{2},\mathcal{B}^1)
 \end{equation}

Let us show that this function is a  $\Gamma$-series in {\it determinants }  $b_i$, $b_{i_1,i_2}$.  The function   $\mathcal{F}$ from  \eqref{bfg} is written as

 \begin{equation}
 (b_{-2}b_{-1})^{\sigma}b_{-1,2}^{2(m_{-2}-k_{-2})}(b_{-2})^{k_{-1}}\cdot \sum_{n\in\mathbb{Z}} \frac{b_{-2,-1}^{\gamma_1+n}}{(\gamma_1+n)!}
 	\frac{b_{-2,1}^{\gamma_2-n}}{(\gamma_2-n)!}
 		\frac{b_{-1}^{2\gamma_3-2n}}{(2\gamma_3-2n)!}
 			\frac{(-b_{1}b_{-1}+b_{-2}b_{2})^{n}}{n!}
 \end{equation}

Consider in detail the sum in this expression. Let us open the brackets in the last factor. Then one gets an expression that is  a  $\Gamma$-series in determinants. It's shift vector  $\delta$ has components

  \begin{align}
 \label{t11}
 \begin{array}{|c|c|c|c|c|c|c|}
 \hline
 &b_{-2,-1} & b_{-2,1} & b_{-1} & -b_{1} & b_{-2} & b_{2}\\
 \hline
\delta & \gamma_1 & \gamma_2& 2\gamma_3 &0 & 0 & 0\\
 \hline
 \end{array}
 \end{align}

And the lattice is generated by two vectors (see row of the table below)



  \begin{align}
 \label{t12}
 \begin{array}{|c|c|c|c|c|c|c|}
 \hline
& b_{-2,-1} & b_{-2,1} & b_{-1} & b_{1} & b_{-2} & b_{2}\\
 \hline
v_1 &1 & -1& -1 &1 & 0 & 0\\
 \hline
 v_0&0 & 0& 1 &1 & -1 & -1\\
 \hline
 \end{array}
 \end{align}
Put

\begin{align}
\begin{split}
\label{dlt}
&\mathcal{B}^{2}=\mathbb{Z}<v_0,v_1>,\\
&\delta=( s_{-2}-m_{-1}  ,   k_{-2}-s_{-2} ,  2(m_{-1}-k_{-1}),0,0,0 )
\end{split}
\end{align}
 One obtains the Lemma.

 \begin{lem}
 \label{lgh}
 	
On the language   $\mathfrak{sp}_4$  a  $h$-highest vector corresponding to a diagram is written as a function

 \begin{equation}
 \label{fb}
 	 (b_{-2}b_{-1})^{\sigma}b_{-1,2}^{k_{-1}-\sigma}(b_{-2})^{2(m_{-2}-k_{-2})}\cdot \mathcal{F}_{\delta}(b_{-2,-1},b_{-2,1},b_{1},b_{-1},b_{-2},b_{2},\mathcal{B}^{2}).
 \end{equation}

 	\end{lem}


\section{The principle Lemma}
\label{oslmma}

In this section we formulate the main instrument that allows us to deal with    $\Gamma$-series.

\begin{lem}
Let us be given a   $\Gamma$-series in determinants $b_{i}$, $b_{i_1,i_2}$ and the underlying lattice has generators  $v_{\alpha}$, $\alpha=1,...,K$,  such that for $\alpha=1,...,k$, $1\leq k\leq K$ the following holds

\begin{enumerate}
\item  
The generators  $v_{\alpha}$ look as follows

$$
v_{\alpha}=e_{X_1}+e_{X_2}-e_{X_3}-e_{X_4},
$$
 where   $e_X$ is a unit vector corresponding to a coordinate  $b_X$,  $X=\{i\}$ or   $\{i_1,i_2\}$.

\item  With each generator $v_{\alpha}$ a relation  between determinants of the following type is associated

 \begin{equation}
 \label{sbv}
b_{X_1}b_{X_2}+b_{X_3}b_{X_4}+b_{X_5}b_{X_6}=0.
\end{equation}

\end{enumerate}

With each generator  $v_{\alpha}$  let us associated a vector

\begin{equation}
\label{ra}
r_{\alpha}=e_{X_1}+e_{X_2}-e_{X_5}-e_{X_6}.
\end{equation}

Then modulo relations   \eqref{sbv} one has

	\begin{equation}
\label{umn}
b_X\mathcal{F}_{\gamma}(b)=\sum_{s\in\mathbb{Z}^k_{\geq 0}} C^{\gamma}_s\mathcal{F}_{\gamma+e_X+sr}(b) \,\,\, mod Pl,
\end{equation}

here  $Pl$ is an ideal generated by  \eqref{sbv}.

\end{lem}

The proof of this Lemma is similar to the proof of an analogous Lemma in the case  $\mathfrak{sp}_4$  from  \cite{a4}  (see also  \cite{a5}).

Let us give explicit formulas for the  coefficients in  \eqref{umn}.  For
  $t,s\in\mathbb{Z}^k_{\geq 0}$ let us introduce functions (we use notations  $tv:=t_1v_1+...+t_kv_k$, $sr:=s_1r_1+...+s_kr_k$ and a multi-index notation for factorials)

\begin{equation}
\mathcal{F}_{\gamma}^s(z)=\sum_{t\in\mathbb{Z}^k}\frac{(t+1)...(t+s)A^{\gamma+tv}}{s!(\gamma+tv)!}
\end{equation}

Than as in  \cite{a4},  one obtains that

	\begin{align}
	\begin{split}
	\label{cs}
	&C^{\gamma}_s=\frac{\mathcal{F}_{\gamma+v}^{s}(1)}{\mathcal{F}^s_{\gamma+v+e_X+sr}(1)}-\sum_{p=0}^{s-1}\frac{\mathcal{F}_{\gamma+v}^{p}(1)
\mathcal{F}_{\gamma+v+pr+e_X+(s-p)r}^{s-p}(1)}{\mathcal{F}_{\gamma+v+pr+e_X+(s-p)r}(1)\mathcal{F}_{\gamma+v+pr+e_X}(1)}=\\
	&=
	\frac{\mathcal{F}_{\gamma+v}^{s}(1)}{\mathcal{F}^s_{\gamma+v+e_X+sr} (1) }-\sum_{p=0}^{s-1}\frac{\mathcal{F}_{\gamma+v}^{p}(1)\mathcal{F}_{\gamma+v+e_X+sr}^{s-p}(1)}{\mathcal{F}_{\gamma+v+e_X+sr}(1)
\mathcal{F}_{\gamma+v+pr+e_X}(1)}
	\end{split}
	\end{align}
\section{A function on  $Sp_4$, corresponding to a diagram}
Let us find a function corresponding to a diagram  \eqref{dc1} on the language   $\mathfrak{sp}_4$.









Note that a  $\Gamma$-series  \eqref{fb}  and the lattice   \eqref{dlt} satisfy the conditions of the principle Lemma and a  relation  \eqref{sbv} corresponds only to   $v_1$.


\subsection{ A function corresponding to a   $h$-highest diagram}
In the space of   $h$-highest vectors there exists a base indexed by  \eqref{dio5},   presented by functions  \eqref{fb}, where   $\delta$  is defined in   \eqref{dlt}.

Let us construct in this space another base also indexed by diagrams.




 To the function   \eqref{fb}  let us apply the ruler  \eqref{umn}, one gets that the function  \eqref{fb} is presented as follows

 \begin{equation}
 const\cdot b_{2,-1}^{m_{-2}} \mathcal{F}_{\omega}(b,\mathcal{B}^2)+l.o.t,
 \end{equation}

where

\begin{equation}
\omega=( s_{-2}-m_{-1}  ,   k_{-2}-s_{-2} ,  2(m_{-1}-k_{-1})+\sigma,0,2k_{-2})-\sigma,0 ),
\end{equation}

and  $l.o.t.$ о is a linear combination of functions of type $\mathcal{F}_{\omega+sr_1}$,  $s\in\mathbb{Z}_{\geq 0}$.


Note that these functions also define   $h$-highest vectors.

Introduce an ordering on the vectors   $\delta$

  $$
  \delta' \prec \delta \Leftrightarrow \delta'=\delta+sr \,\,\,mod \,\,\,\mathcal{B}^2
  $$

 Then one gets that  in the space of    $h$-highest vectors one has two collections of functions indexed by diagrams: $\{b_{2,-1}^{m_{-2}-k_{-2}} \mathcal{F}_{\omega}(b,\mathcal{B}^{2}) \}$,  and functions  \eqref{fb}. The functions from the second collection are expressed from function of the first collection using a triangular operator. The second collection is a base in the space of    $h$-highest vectors  thus the first collection is also a base in this space.

Below   in the space of $h$-highest vectors with a fixed   $h$-weight   the base $b_{2,-1}^{m_{-2}-k_{-2}} \mathcal{F}_{\omega}(b,\mathcal{B}^{2})$ is considered.

  \subsection{ A pseudo-action of  $\mathfrak{gl}_4$}

Determinants   $b_{i}$, $b_{i_1,i_2}$ satisfy some relations. These are the Plucker relations and also for example  the relation  $b_{-1,1}=-b_{-2,2}$.  Consider the symbols  $ B_{i}$, $B_{i_1,i_2}$,  that satisfy   the Plucker relations   {\it only  }.  Also they are supposed to be antisymmetric  under permutations of  $i_1$ and   $i_2$.

 These symbols ban be identified with function on  $GL_4$.

One has an action of  $\mathfrak{gl}_4$ onto these symbols

 $$
E_{i,j}B_{i_1}=B_{i_1\mid_{j\mapsto i }}.
 $$

These formulas do not give an action of   $\mathfrak{gl}_4$ onto determinants   $b_{i}$, $b_{i_1,i_2}$.  But the formula
	
	 $$f_{i,j}=E_{i,j}+sign(i)sign(j) E_{-j,-i}$$ gives an action of   $\mathfrak{sp}_4$ onto  $B$.

\subsection{A function corresponding to a diagram}

Consider an expression
 
 \begin{equation}
 \label{fb2}
B_{2,-1}^{m_{-2}-k_{-2}} \mathcal{F}_{\omega}(B,\mathcal{B}^{2}).
 \end{equation}
 

Apply the operator   \begin{equation}\label{ff}\frac{f^p_{1,-2}}{p!}=\sum_{p_1+p_2=p}
\frac{E_{2,-1}^{p_1}}{p_1!}\frac{E_{1,-2}^{p_2}}{p_2!}.\end{equation}

Here we use the fact that   $E_{1,-2}$ and  $E_{2,-1}$  commute and they do not act onto   $B_{2,-1}^{m_{-2}-k_{-2}}$.

\begin{remark}
	 It one tries to apply   $\frac{F_{0,-2}^p}{p!}$ to the function  \eqref{bbo},  that it is not possible to write an analogous expression since the operators   $E_{0,-2}$ and   $E_{2,0}$ do not commute.   This is the main  obstacle in the case of the language  $\mathfrak{o}_5$.
\end{remark}

Let us write the operator  $E_{1,-2}$ explicitly

\begin{equation}
\label{oe1}
B_{1}\frac{\partial }{\partial B_{-2}}+B_{1,-1}\frac{\partial }{\partial B_{-2,-1}}+B_{1,2}\frac{\partial }{\partial B_{-2,2}}.
\end{equation}

Analogously one writes   $E_{2,-1}$ as follows:

\begin{equation}
\label{oe2}
B_{2}\frac{\partial }{\partial B_{-1}}+B_{-2,2}\frac{\partial }{\partial B_{-2,-1}}+B_{1,2}\frac{\partial }{\partial B_{1,-1}}.
\end{equation}

The summand in the expressions    \eqref{oe1} and  \eqref{oe2} commute.

Instead of the symbol  $B_1$  let us introduce symbols  $B'_1$ and  $B''_1$. In the function $\mathcal{F}_{\delta}$,  corresponding to a  $h$-highest vector let us change  $B_1$ to   $B'_1$. And in the operator  \eqref{oe1} let us change  $B_1$ to  $B''_1$.

Analogously one introduces   $B'_2$ and   $B''_2$.  One writes  $B'_2$ instead of   $B_2$ in  $\mathcal{F}_{\delta}$, and one writes  $B''_2$ instead of  $B_2$  in  \eqref{oe2}.

The operator   $\frac{\partial }{\partial B_{X}}$ acts onto a   $\Gamma$-series as a transformation  $\omega \mapsto \omega-e_{X}$.
Then we can write:

\begin{align}
\begin{split}
& \frac{f_{1,-2}^p}{p!}\mathcal{F}_{\omega}
=\sum_{p_1+p_2+q_1+q_2+q_3=p}
\frac{(B''_1)^{p_1}}{p_1!}
\frac{B_{1,-1}^{p_2-q_3}}{(p_2-q_3!)}\frac{(B''_2)^{q_1}}{q_1!}
\frac{B_{-2,2}^{q_2}}{q_2!}\frac{B_{1,2}^{q_3}}{q_3!}\cdot\\
 & \cdot \mathcal{F}_{\omega-p_1e_{-2}-(p_2+q_2)e_{-2,-1}-q_1e_{-1}}
(B_{-2,-1},B_{-2,1},B_{-1},-B'_{1},B_{-2},B'_{2},B^2)
\end{split}
\end{align}

This expression is a   $\Gamma$-series in variables  $$ B''_1,B'_2, B''_2,B_{-1,1},B_{-2,2},B_{1,2},B_{-2,-1},B_{-2,1}, B_{-1}, -B'_1, B_{-2}, B'_{1}.$$

The shift vector looks as follows

\begin{equation}
(0,0,0,p,0,0, \delta-pe_{-2,-1})
\end{equation}

Let us write the generators of the lattice. These are the vectors  involved in construction of   $ \mathcal{F}_{\omega}$ and also the vectors defining the shifts of exponents corresponding to the change in  $p_1,...,q_3$  preserving the sum  $p_1+p_2+p_3+q_1+q_2+q_3=p$.  Such changes are generated by

\begin{equation}
\begin{array}{|c|c|c|c|c|c|}
\hline
 & p_1 & p_2 & q_1 & q_2 & q_3\\
\hline
u_0 & 0 &1 &0 &-1&0\\
\hline
u_1 & 0 &0 &0 &-1&1\\
\hline
u_2 & 0 &0 &1 &0&-1\\
\hline
u_3 & 1 &-1 &0 &0&0\\
\hline
\end{array}
\end{equation}

Using these vectors let us  obtain vectors of changes of exponents of  $B_X$.  For the vector  $u_0$  one obtains a vector

\begin{equation}
w_0=e_{-1,1}-e_{-2,2},
\end{equation}

for vectors   $u_i$, $i=1,2,3$  one obtais vectors:

 \begin{align}
 	\begin{split}
 		& v_2=e_{-2,-1}+e_{1,2}-e_{-1,1}-e_{-2,2},\\
 		& v_3=e''_{2}+e_{-1,1}-e_{-1}-e_{2,1},\\
 		& v_4=e''_{1}+e_{-2,-1}-e_{-2}-e_{-1,1}.
 	\end{split}
\end{align}

 Thus one obtains the Theorem.

 \begin{thm}
 In the space of an irreducible representation of    $g$  there exists a Gelfand-Tsetlin type base corresponding to the chain  $g\supset h$, indexed by diagrams  \eqref{dc1},  to a diagram there corresponds the function
 	\begin{align}
 	\begin{split}
 \label{omg}
 	& b_{2,-1}^{m_{-2}-k_{-2}}\mathcal{F}_{\omega}(b,\mathcal{B}_{GC})\mid_{b'_1=b''_1=b_1,\,\,\, b'_2=b''_2=b_2},\\
 	& b=(b''_1, b''_2,b_{-1,1},b_{-2,2},b_{1,2},b_{-2,-1},b_{-2,1}, b_{-1}, -b'_1,b_{-2}, b'_{2})	,\\
&\omega=(0,0,s_{-2}-s_{-1},0,0,s_{-1}-m_{-1}  ,   k_{-2}-s_{-2} ,  2(m_{-1}-k_{-1})+\sigma,0,2k_{-2}-\sigma,0 ) ,\\
& \mathcal{B}_{GC}=\mathbb{Z}	<v_0,w_0,v_1,v_2,v_3,v_4>.
 	\end{split}
 	\end{align}
 \end{thm}

\begin{remark} The vector  $\omega$ satisfies the following selection rulers:

$$
\omega=(0,0,\omega_{-1,1},0,0,\omega_{-2,-1},\omega_{-2,1},\omega_{-1},0,\omega_{-2},0).
$$

\end{remark}

\section{The action of generators}

\subsection{The action of an arbitrary $f_{i,j}$}

Consider an arbitrary generator $f_{i,j}$ of the algebra  $\mathfrak{sp}_4$. Apply it to  $\mathcal{F}_{\omega}$, one gets

\begin{align*}
& f_{i,j}\mathcal{F}_{\omega}=\sum_{Y}b_{i,Y}\frac{\partial}{\partial b_{j,Y}}\mathcal{F}_{\omega}-sign(i)sign(j)\sum_{Z}b_{-j,Z}\frac{\partial}{\partial b_{-i,Z}}\mathcal{F}_{\omega}
\end{align*}

A summation is taken over all subsets   $Y,X\subset\{\pm 1,\pm 2\}$, such that  $i,j\notin Y$, $-i,-j\notin Z$.  A differentiation acts onto a    $\Gamma$-series by a change of  a  shift-vector, and a multiplication by a determinants acts according to   \eqref{umn}. One has

\begin{align}
\begin{split}
\label{fijd}
&f_{i,j}\mathcal{F}_{\omega}=\sum_{Y}\sum_{s^1\in\mathbb{Z}^{k}_{\geq 0}}C_{s^1}^{\omega-e_{Y,j}}\mathcal{F}_{\omega+e_{Y,i}-e_{Y,j}+s^1r}-\\&-sign(i)sign(j)
\sum_{Z}\sum_{s^2\in\mathbb{Z}^{k}_{\geq 0}}C_{s^2}^{\omega-e_{Z,-i}}\mathcal{F}_{\omega+e_{Z,-j}-e_{Z,-i}+s^2r},
\end{split}
\end{align}

where $C_{s}^{\omega}$ is given by   \eqref{cs}.

A transformation of the shift-vector $$\omega \mapsto \omega+e_{Y,i}-e_{Y,j}+s^1r,\,\,\, \omega\mapsto \omega+e_{Z,-j}-e_{Z,-i}+s^2r$$
corresponds to a transformation of the Gelfand-Tsetlin  diagram.  Let us derive it for some generators of   $\mathfrak{sp}_4$.
It is  sufficient to do it for  $\frac{1}{2}(f_{-2,-2}+f_{-1,-1})$ and   $\frac{1}{2}(f_{-2,-2}-f_{-1,-1})$, $f_{1,-2}$ and   $f_{-2,1}$, $f_{-2,-1}$ and  $f_{-1,-2}$.

\subsection{The action of Cartan elements}
The Cartan elements $\frac{1}{2}(f_{-2,-2}+f_{-1,-1})$ and  $\frac{1}{2}(f_{-2,-2}-f_{-1,-1})$ on the language of   $\mathfrak{o}_5$ turn into  $F_{-2,-2}$ and  $F_{-1,-1}$.

 The Gelfand-Tsetlin type base consists of eigenvectors of   $F_{-1,-1}$, $F_{-2,-2}$ and the eigenvalues equal to

\begin{align*}
&\lambda_{-2}=\sigma+2(k_{-2}+k_{-1})-(m_{-2}+m_{-1})-s_{-2},\\
&\lambda_{-1}=s_{-1}.
\end{align*}

This fact was derived in  \ref{o335}.

 \subsection{The action of $f_{1,-2}$ and  $f_{-2,1}$ }
The function corresponding to a diagram is obtained from a  function corresponding to a   $h$-highest diagram by the action of   $\frac{f_{1,-2}^{s_{-2}-s_{-1}}}{(s_{-2}-s_{-1})!}$

Hence the action of   $f_{1,-2}$  gives a transformation   $s_{-1}\mapsto s_{-1}-1$ and a multiplication onto  $s_{-2}-s_{-1}+1$.

 The action of  $f_{-2,1}$  gives a transformation  $s_{-1}\mapsto s_{-1}+1$  and a multiplication onto   $s_{-2}-s_{-1}+1$.

\subsection{The action of $f_{-1,1}$ and   $f_{-2,2}$}

Let us write  $f_{-1,1}$  as

$$
b_{-1}\frac{\partial}{\partial b'_{1}}+b_{-1}\frac{\partial}{\partial b''_{1}}+
b_{-2,-1}\frac{\partial}{\partial b_{-2,1}}+b_{-1,2}\frac{\partial}{\partial b_{1,2}}
$$

Also write  $f_{-2,2}$ as

$$
b_{-2}\frac{\partial}{\partial b'_{2}}+b_{-2}\frac{\partial}{\partial b''_{2}}+b_{-2,-1}\frac{\partial}{\partial b_{2,- 1}}+
b_{-2,1}\frac{\partial}{\partial b_{2,1}}+
$$

The action of these expressions onto   $\mathcal{F}_{\delta}$ is given by  \eqref{fijd}.
Let us find a transformation of a diagram  corresponding to the  transformation  $\omega\mapsto \omega-e_{X}+e_{Y}$ of the shift vector associated with the  summand  $b_{Y}\frac{\partial}{\partial b_{X}}$. Let us also describe a transformation of a diagram corresponding to the addition of  $r_{\alpha}$.


\subsubsection{  The addition of $r_{\alpha}$}

With each vector   $v_{\alpha}$ one associates a  vector  $r_{\alpha}$,  defined in   \eqref{ra}.
Consider the vector   $$v_1=e_{-2,-1}-e_{-2,1}-e_{-1}+e_{1}\,\,\Rightarrow  \,\, r_{1}=e_{-2}+e_{-1,1}-e_{-2,1}-e_{-1}.$$

Addition  of    $r_{1}$ to  $\omega$ from  \eqref{omg} is equivalent   $mod\mathcal{B}_{GC}$ to the following transformation:

$$
\begin{cases}
s_{-1}\mapsto s_{-1}+1,\\
\sigma=0,\,\,\,k_{-1}\mapsto \sigma=1,\,\,\,k_{-1}+1\\
\sigma=1,\,\,\,k_{-1}\mapsto \sigma=0,\,\,\,k_{-1}
\end{cases}
$$


Consider the vector  $$v_{2}=e_{-2,-1}+e_{1,2}-e_{-1,1}-e_{-2,2}\,\,\Rightarrow  \,\,
r_{2}=e_{-2,1}+e_{1,-2}-e_{-1,1}-e_{-2,2}
$$

The addition of    $r_{2}$ to  $\omega$ from  \eqref{omg} is equivalent   $mod\mathcal{B}_{GC}$ to the following transformation:
$$
\begin{cases}
k_{-2}\mapsto k_{-2}-1,\,\,
s_{-2}\mapsto s_{-2}-2
\end{cases}
$$


Consider the vector  $$v_{3}=e''_{2}+e_{-1,1}-e_{-1}-e_{1,2}\,\,\Rightarrow  \,\,
r_{3}=e''_{1}+e_{-1,2}-e_{-1}-e_{1,2}
$$

The addition of    $r_{3}$ to  $\omega$ from  \eqref{omg} is equivalent   $mod\mathcal{B}_{GC}$ to the following transformation:

$$
\begin{cases}
k_{-2}\mapsto k_{-2}-1,\\
s_{-2}\mapsto s_{-2}-1,\\
\sigma=0,\,\,\,k_{-1}\mapsto \sigma=1,\,\,\,k_{-1}+1\\
\sigma=1,\,\,\,k_{-1}\mapsto \sigma=0,\,\,\,k_{-1}
\end{cases}
$$








Consider the vector   $$v_5=e''_{1}+e_{-2,-1}-e_{-1,1}-e_{-2}\,\,\Rightarrow  \,\,
r_{5}=e_{-1}+e_{-2,1}-e_{-1,1}-e_{-2}
$$

The addition of    $r_{4}$ to  $\omega$ from  \eqref{omg} is equivalent   $mod\mathcal{B}_{GC}$ to the following transformation:

$$
\begin{cases}
s_{-2}\mapsto s_{-2}-1,\\
\sigma=1,\,\,\,k_{-1}\mapsto \sigma=0,\,\,\,k_{-1}-1\\
\sigma=0,\,\,\,k_{-1}\mapsto \sigma=1,\,\,\,k_{-1}
\end{cases}
$$

\subsubsection{ The transformation  $\omega-e_{Y,1}+e_{Y,-1}$  associated with the action of  $f_{-1,1}$}

Consider the summand  $b_{-1}\frac{\partial}{\partial b'_{1}}$.  The  transformation of the shift vector  looks as follows:   $\omega\mapsto \omega+e_{-1}-e'_{1}$.
 It is equivalent   $mod\mathcal{B}_{GC}$ to the following transformation of the diagram:

$$
 \begin{cases} s_{-2}\mapsto s_{-2}+1,\,\, s_{-1}\mapsto s_{-1}+1\end{cases}
$$

Consider the summand  $b_{-1}\frac{\partial}{\partial b''_{1}}$.  The  transformation of the shift vector  looks as follows:   $\omega\mapsto \omega+e_{-1}-e''_{1}$.
 It is equivalent   $mod\mathcal{B}_{GC}$ to the following transformation of the diagram:

$$
\begin{cases} s_{-2}\mapsto s_{-2}+1 \\
\sigma=1,\,\,\,k_{-1}\mapsto \sigma=0,\,\,\,k_{-1}-1\\
\sigma=0,\,\,\,k_{-1}\mapsto \sigma=1,\,\,\,k_{-1}
\end{cases}
$$

Consider the summand   $b_{-2,-1}\frac{\partial}{\partial b_{-2,1}}$. The  transformation of the shift vector  looks as follows:    $\omega\mapsto \omega+e_{-2,-1}-e_{-2,1}$.
 It is equivalent   $mod\mathcal{B}_{GC}$ to the following transformation  of the diagram:

$$
\begin{cases} s_{-2}\mapsto s_{-2}+1,\,\,s_{-1}\mapsto s_{-1}+1
\end{cases}
$$

Consider the summand  $b_{-1,2}\frac{\partial}{\partial b_{1,2}}$. The  transformation of the shift vector  looks as follows:   $\omega\mapsto \omega+e_{-1,2}-e_{1,2}$.
 It is equivalent   $mod\mathcal{B}_{GC}$ to the following transformation  of the diagram:

$$
\begin{cases}
k_{-2}\mapsto k_{-2}+1,\,\,
 s_{-2}\mapsto s_{-2}+1 \,\,s_{-1}\mapsto s_{-1}-1
\end{cases}
$$

\subsubsection{ The transformation  $\omega-e_{Y,2}+e_{Y,-2}$  associated with the action of $f_{-2,2}$}

Consider the summand $b_{-2}\frac{\partial}{\partial b'_{2}}$. The  transformation of the shift vector  looks as follows:   $\omega\mapsto \omega+e_{-2}-e'_{2}$.
 It is equivalent   $mod\mathcal{B}_{GC}$ to the following transformation  of the diagram:

$$
\begin{cases}
k_{-1}\mapsto k_{-1}+1,\,\,
s_{-2}\mapsto s_{-2}+1 \,\,s_{-1}\mapsto s_{-1}+1
\end{cases}
$$

Consider the summand  $b_{-2}\frac{\partial}{\partial b''_{2}}$.  The  transformation of the shift vector  looks as follows:    $\omega\mapsto \omega+e_{-2}-e''_{2}$.
 It is equivalent   $mod\mathcal{B}_{GC}$ to the following transformation  of the diagram:

$$
\begin{cases}
s_{-1}\mapsto s_{-1}-1,\\
\sigma=0,\,\,\,k_{-1}\mapsto \sigma=1,\,\,\,k_{-1}+1\\
\sigma=1,\,\,\,k_{-1}\mapsto \sigma=0,\,\,\,k_{-1}
\end{cases}
$$

Consider the summand $b_{-2,-1}\frac{\partial}{\partial b_{2,-1}}$.  The  transformation of the shift vector  looks as follows:   $\omega\mapsto \omega+e_{-2,-1}-e_{2,-1}$.
 It is equivalent   $mod\mathcal{B}_{GC}$ to the following transformation  of the diagram:

$$
\begin{cases}
k_{-2}\mapsto k_{-2}+1,\,\,
s_{-2}\mapsto s_{-2}+1,\,\,
s_{-1}\mapsto s_{-1}+1
\end{cases}
$$

Consider the summand $b_{-2,1}\frac{\partial}{\partial b_{2,1}}$.  The  transformation of the shift vector  looks as follows:    $\omega\mapsto \omega+e_{-2,1}-e_{2,1}$.
 It is equivalent   $mod\mathcal{B}_{GC}$ to the following transformation of the diagram:

$$
\begin{cases}
s_{-2}\mapsto s_{-2}+2,\,\,s_{-1}\mapsto s_{-1}+2
\end{cases}
$$

\end{document}